\documentclass{article}

\def\ra{\rightarrow}

\def\ss{\subseteq}

\def\Aut{{\rm Aut}}
\def\boac{{\rm boundary orbit accumulation point}}

 \def\HollowBox #1#2{{\dimen0=#1 \advance\dimen0 by -#2       
       \dimen1=#1 \advance\dimen1 by #2                       
        \vrule height #1 depth #2 width #2                    
        \vrule height 0pt depth #2 width #1                   
        \llap{\vrule height #1 depth -\dimen0 width \dimen1}%
       \hskip -#2                                             
       \vrule height #1 depth #2 width #2}}                   
 \def\BoxOpTwo{\mathord{\HollowBox{6pt}{.4pt}}\;}             

\def\endpf{\hfill $\BoxOpTwo$}

\font\teneufm=eufm10
\font\seveneufm=eufm7
\font\fiveeufm=eufm5
\newfam\eufmfam
\textfont\eufmfam=\teneufm
\scriptfont\eufmfam=\seveneufm
\scriptscriptfont\eufmfam=\fiveeufm

\newfam\msbfam
\font\tenmsb=msbm10  \textfont\msbfam=\tenmsb
\font\sevenmsb=msbm7  \scriptfont\msbfam=\sevenmsb
\font\fivemsb=msbm5    \scriptscriptfont\msbfam=\fivemsb
\def\Bbb{\fam\msbfam \tenmsb}

\def\CC{{\Bbb C}}

\newtheorem{theorem}{Theorem}[section]

\newtheorem{proposition}[theorem]{Proposition}
\newtheorem{lemma}[theorem]{Lemma}
\newtheorem{remark}[theorem]{Remark}

\makeindex

\usepackage{graphicx}

\usepackage{amsmath}

\begin{document}

\begin{center}
\huge \bf
The Automorphism Group of a Domain with an Exponentially Flat Boundary Point\footnote{{\bf Key Words:}  orbit accumulation point,
automorphism, orbit, holomorphic mapping.}\footnote{{\bf MR Classification
Numbers:}  32M05,32M99.}
\end{center}
\vspace*{.12in}

\begin{center}
\large Steven G. Krantz\footnote{Author supported in part
by the National Science Foundation and by the Dean of the Graduate
School at Washington University.}\end{center}
\vspace*{.15in}

\begin{center}
\today
\end{center}
\vspace*{.2in}

\begin{quotation}
{\bf Abstract:} \sl
We study the automorphisms group action on a bounded
domain in $\CC^n$ having a boundary point that is exponentially
flat.   Such a domain typically has a compact automorphism group.
Our results enable us to generate many new examples.
\end{quotation}
\vspace*{.25in}

\setcounter{section}{-1}

\section{Introduction}

A {\it domain} $\Omega$ in $\CC^n$ is a connected, open set.  An {\it automorphism}
of $\Omega$ is a biholomorphic self-map.  The collection of automorphisms forms
a group under the binary operation of composition of mappings.   The topology
on this group is uniform convergence on compact sets, or the compact-open topology.
We denote the automorphism group by $\Aut(\Omega)$.

Although domains with {\it transitive automorphism group} are of some interest,
they are relatively rare (see [HEL]).   A geometrically more natural condition
to consider, and one that gives rise to a more robust and broader class of domains,
is that of having {\it non-compact automorphism group}.   Clearly a domain has
non-compact automorphism group if there are automorphisms $\{\varphi_j\}$ which
have no subsequence that converges to an automorphism.   The following proposition
of Henri Cartan is of particular utility in the study of these domains:

\begin{proposition} \sl
Let $\Omega \ss \CC^n$ be a bounded domain.  Then $\Omega$ has non-compact automorphism
group if and only if there are a point $X \in \Omega$, a point $P \in \partial \Omega$, and
automorphisms $\varphi_j$ of $\Omega$ such that $\varphi_j(X) \ra P$ as $j \ra \infty$.
\end{proposition}

We refer the reader to [NAR] for discussion and proof of Cartan's result.

A point $P$ in $\partial \Omega$ is called a {\it boundary orbit accumulation point} if there is an $X \in \Omega$ and
automorphisms $\varphi_j$ of $\Omega$ such that $\lim_{j\ra \infty} \varphi_j(X) = P$.

In the paper [GRK1], we considered the domain
$$
\Omega_\infty = \{(z_1, z_2) \in \CC^2: |z_1|^2 + 2 \exp(-1/|z_2|^2) < 1\} \, .     \eqno (*)
$$
We showed that this domain has {\it compact automorphism group}.  In particular,
the only automorphisms (i.e., biholomorphic self-maps) of $\Omega$ are the rotations
in each variable separately.  We note that, unlike the domains
$$
\Omega_m \equiv \{(z_1, z_2); |z_1|^2 + |z_2|^{2m} < 1\}, m = 1, 2, \dots,
$$
which have all points of the form $(e^{i\theta}, 0)$ as boundary orbit
accumulation points, $\Omega_\infty$ has {\it no} boundary orbit
accumulation points.

An example such as this bears directly on the Greene-Krantz conjecture:

\begin{quote}
{\bf Conjecture:}  Let $\Omega \ss \CC^n$ be a smoothly bounded domain.  If $P$ is
a \boac, then $P$ must be a point
of finite type in the sense of Kohn/Catlin/D'Angelo (see [KRA1] for this
concept).
\end{quote}

This conjecture is not known to be true in its full generality.   But evidence for
its correctness is provided, for instance, in [KIM] and [KIK].

The domain $\Omega_\infty$ in $(*)$ can be analyzed like this:  If there is
an $X \in \Omega_\infty$ and a $P = (p_1, p_2) \in \partial \Omega_\infty$ with $p_2 \ne 0$
so that $\lim_{j \ra \infty} \varphi_j(X) = P$ for some automorphisms $\varphi_j$,
then $P$ is strongly pseudoconvex.  The theorem of Bun Wong and Rosay (see [WON], [ROS])
then tells us that $\Omega_\infty$ must be biholomorphic to the ball.  But a result of
Bell/Boas [BEB] tells us that such a biholomorphism must extend smoothly to the boundary.
That is impossible, since $\Omega_\infty$ has a circle $\{(e^{i\theta},0)\}$ of weakly
pseudoconvex points in the boundary while the ball $B$ is strongly pseudoconvex.
If instead $p_2 = 0$, then a delicate analysis---see the proof of Proposition 1.3 
below---shows that it is impossible for
an orbit $\{\varphi_j(X)\}$ to accumulate at $P$.  Thus we conclude, using
Proposition 0.1, that $\Aut(\Omega_\infty)$ must be compact.

Of course the boundary points $(e^{i\theta}, 0)$ of the domain $\Omega_\infty$
in $(*)$ are of infinite type.   It is important to have examples like
this at hand to aid in the study of the Greene-Krantz conjecture.

The purpose of this paper is to provide a good many further examples of such domains.
We provide domains with different types of geometry in different dimensions.  All have
compact automorphisms group.

In general it is quite difficult to produce examples of any kind in the study
of automorphism groups of domains in $\CC^n$.  The techniques presented in this
paper may prove useful in other contexts.

\section{Some Principal Results}

The following simple technical result will be useful.  I thank John P. D'Angelo for
helpful discussions of the idea.

\begin{lemma} \sl
Let $F$ be a holomorphic function on the unit ball $B$ in $\CC^n$, $n \geq 2$, which
extends continuously to $\overline{B}$.	 Assume that $|F|$ on $\partial B$ equals
a positive constant $c$.   Then $F$ is identically constant.
\end{lemma}
{\bf Proof:}   Let $P$ be a boundary point of the ball.  Let {\bf d} be an analytic
disc whose image lies in a complex line which intersects $\partial B$ transversally.
Also suppose that the
boundary circle of this analytic disc lies in $\partial B$---in other words, the analytic disc
is just a complex line intersected with $\overline{B}$.  Assume that {\bf d} is very close to $P$.
Then $|F|$ restricted to {\bf d} will assume values that are uniformly very close to c.
But $|F|$ restricted to {\bf d} is a function of one complex variable that satisfies the
hypotheses of the lemma.  So $F$ restricted to {\bf d} must be a finite Blaschke product
with values which in modulus are very close to $c$.  It follows that $F$ can have no zeros,
so $F$ restricted to {\bf d} must be a constant.   Since this statement must hold on all
such analytic discs closed to $P$, we may conclude that $F$ is constant.
\endpf
\smallskip \\

\begin{remark} \rm
In point of fact, this lemma is true on any smoothly bounded domain.   For such a domain
will alway have a relatively open boundary neighborhood that is strongly convex.  And the
argument just presented will be valid on that neighborhood.
\end{remark}

\begin{proposition} \sl
Consider the domain
$$
\widetilde{\Omega} = \{z = (z_1, z_2, \dots, z_n) \in \CC^n: |z_1|^2 + |z_2|^2 + \cdots + |z_{n-1}|^2 + \psi(|z_n|) < 1\} \, ,
$$
where $\psi$ is a real-valued, even, smooth, monotone-and-convex-on-$[0,\infty)$ function of a real variable  with $\psi(0) = 0$ 
that vanishes to infinite order at 0.   Then there
do {\it not} exist a point $X \in \widetilde{\Omega}$ and automorphisms $\varphi_j$ of $\widetilde{\Omega}$ so that
$\lim_{j \ra \infty} \varphi_j(X) = (1, 0, \dots, 0) \in \partial \widetilde{\Omega}$.
\end{proposition}

\begin{remark} \rm
Observe that points $\zeta$ of $\partial \widetilde{\Omega}$ with $\zeta_n \ne 0$ are strongly pseudoconvex (as a calculation
shows).  Such a point cannot be an orbit accumulation point because then, by the Bun Wong/Rosay theorem, the
domain would have to be biholomorpic to the unit ball.

By contrast, the boundary points of the form $(e^{i\theta},0,\dots,0)$, $(0, e^{i\theta}, 0, \dots, 0)$,
\dots, $(0,0, \dots, e^{i\theta},0)$ are of infinite type.   We state
a result for the particular boundary point $(1,0,0, \dots, 0)$, but in fact the analysis applies to
any boundary point of the form $(e^{i\theta}, 0,0,\dots, 0)$, $(0, e^{i\theta}, 0, \dots, 0)$,
\dots, $(0,0, \dots, e^{i\theta},0)$,$0 \leq \theta \leq 2\pi$.

Of course an example of a real function $\psi$ as in the statement of the theorem
is $\lambda(t) = 2\exp(-1/t^2)$.
\end{remark}

\noindent {\bf Proof of Proposition 1.3:}   Seeking a contradiction, let us suppose
that there is a point $X \in \widetilde{\Omega}$ and a sequence of automorphisms $\varphi_j$ such
that $\varphi_j(X) \ra (1,0,\dots, 0) \in \partial \widetilde{\Omega}$.

The set of weakly pseudoconvex points in the boundary
is 
$$
S = \{z \in \partial \widetilde{\Omega}: z_n = 0\} \, .
$$
All other boundary points are strongly pseudoconvex.  Since the work of Bell/Boas [BEB] tells us
that automorphisms extend smoothly to the boundary, we may conclude that any automorphism $\varphi$ will preserve
$S$.   Thus it will also preserve the ball
$$
{\bf e} = \{(z_1, z_2, \dots, z_{n-1}, 0): \sum_{j=1}^{n-1} |z_j|^2 \leq 1\} \, .
$$

We see then that the restriction of $\varphi$ to {\bf e} will be an automorphism of
the unit ball in $\CC^{n-1}$.  After composing each $\varphi_j$ with rotations in the
$z_1$, $z_2$, \dots, $z_{n-1}$ variables, we may take it that
\begin{eqnarray*}
\varphi_j(z_1, z_2, z_3\dots, z_{n-1} ,0) & = & \hbox{ \ \  \ \ \ \ \ \ \ \ \ \ \ \ \ \ \ \ \ \ \ \ \ \ \ \ \ \ \ \ \ \ \ \ \ \ \ \ \ \ \ \ \ \ \ \ \ \ } \\
\end{eqnarray*}
\vspace*{-.3in}
$$
\left ( \frac{z_1 - a^1_j}{1 - \overline{a}^1_j z_1}, \frac{\sqrt{1 - |a^2_j|^2}z_2}{1 - \overline{a}^2_j z_1}, \frac{\sqrt{1 - |a^3_j|^2}z_3}{1 - \overline{a}^3_j z_1}, \dots
     \frac{\sqrt{1 - |a^{n-1}_j|^2}z_{n-1}}{1 - \overline{a}^{n-1}_j z_1}, 0 \right ) \, .
$$
for some $a_j^\ell \in \CC$, $|a_j^\ell| < 1$, $\ell = 1, \dots, n-1$.  By imitating the proof of the classical result (see [KRA1]) about automorphisms of
circular domains, we may see that each $\varphi_j$ commutes with rotations in the $z_n$ variable.  It follows
that each $\varphi_j$ preserves every disc of the form
$$
{\bf f}_\alpha \equiv \{(\alpha, \zeta): |\alpha|^2 + \psi(|\zeta|) < 1\} \, ,
$$
for some fixed $\alpha \in \CC^{n-1}$ with $|\alpha| < 1$.  Of course here $\zeta \in \CC$. By rotational invariance, the image
of ${\bf f}_\alpha$ under $\varphi_j$ must be a disc of the form
$$
\varphi_j({\bf f}_\alpha) = \{(\varphi_j(\alpha,0), \zeta): |\varphi_j(\alpha,0)|^2 + \psi(|\zeta|) < 1\} \, .
$$
Furthermore, $\varphi_j$ sends the center of ${\bf f}_\alpha$ to the center of the image disc.  We conclude
then that $\varphi_j$ must have the form
\begin{eqnarray*}
(z_1, z_2, z_3, \dots, z_n) & \longmapsto & \hbox{ \ \  \ \ \ \ \ \ \ \ \ \ \ \ \ \ \ \ \ \ \ \ \ \ \ \ \ \ \ \ \ \ \ \ \ \ \ \ \ \ \ \ \ \ \ \ \ \ } \\
\end{eqnarray*}
\vspace*{-.3in}
$$
\left ( \frac{z_1 - a^1_j}{1 - \overline{a}^1_j z_1}, \frac{\sqrt{1 - |a^2_j|^2}z_2}{1 - \overline{a}^2_j z_1}, \frac{\sqrt{1 - |a^3_j|^2}z_3}{1 - \overline{a}^3_j z_1}, \dots
     \frac{\sqrt{1 - |a^{n-1}_j|^2}z_{n-1}}{1 - \overline{a}^{n-1}_j z_1}, z_n \cdot \lambda_j(z_1, \dots,z_{n-1}) \right ) \, ,
$$
for $\lambda_j$ holomorphic.

Thus we see that
\begin{align}
\left | \frac{z_1 - a_j^1}{1 - \overline{a}_j^1 z_1} \right |^2 & +  \left | \frac{\sqrt{1 - |a_j^2|^2} z_2}{1 - \overline{a}_j^2 z_1} \right |^2
 + \left | \frac{\sqrt{1 - |a_j^3|^2} z_3}{1 - \overline{a}_j^3 z_1} \right |^2 + \cdots +  \notag \\
& \qquad \left | \frac{\sqrt{1 - |a_j^{n-1}|^{n-1}} z_2}{1 - \overline{a}_j^{n-1} z_1} \right |^2 + \psi(z_n \cdot\lambda_j(z_1,\dots,z_{n-1})) < 1 \, .   \tag{$\dagger$}  \\ \notag
\end{align}
This may be rewritten as
\begin{eqnarray*}
|z_n| < \frac{1}{|\lambda_j(z_1, z_2, \dots, z_{n-1})|} \cdot \psi^{-1} \Biggl ( 1 & - & 
\left | \frac{z_1 - a_j^1}{1 - \overline{a}_j^1 z_1} \right |^2 \\
 \qquad - \left | \frac{\sqrt{1 - |a_j^2|^2} z_2}{1 - \overline{a}_j^2 z_1} \right |^2
 & - & \left | \frac{\sqrt{1 - |a_j^3|^2} z_3}{1 - \overline{a}_j^3 z_1} \right |^2 - \cdots -   
   \left | \frac{\sqrt{1 - |a_j^{n-1}|^2} z_{n-1}}{1 - \overline{a}_j^{n-1} z_1} \right |^2 \Biggr ) \\
\end{eqnarray*}

But we also know that
$$
|z_n| < \psi^{-1} \left ( 1 - |z_1|^2 - |z_2|^2 - \cdots - |z_{n-1}|^2 \right ) \, .
$$
The only possible conclusion is that
$$
\frac{\psi^{-1} \Biggl ( 1  -  
\left | \frac{z_1 - a_j^1}{1 - \overline{a}_j^1 z_1} \right |^2 
  - \left | \frac{\sqrt{1 - |a_j^2|^2} z_2}{1 - \overline{a}_j^2 z_1} \right |^2
  -  \left | \frac{\sqrt{1 - |a_j^3|^2} z_3}{1 - \overline{a}_j^3 z_1} \right |^2 - \cdots -   
   \left | \frac{\sqrt{1 - |a_j^{n-1}|^2} z_{n-1}}{1 - \overline{a}_j^{n-1} z_1} \right |^2 \Biggr )}{\psi^{-1} \left ( 1 - |z_1|^2 - |z_2|^2 - \cdots - |z_{n-1}|^2 \right )} 
$$
$$
 \hbox{ \ \ } \qquad \qquad \qquad \qquad \qquad \qquad \qquad  = |\lambda_j(z_1, z_2, \dots, z_{n-1})| \, .  \eqno (\star)
$$
Now we let $|(z_1, z_2, \dots, z_{n-1})| \ra 1$. Since $\psi^{-1}$ vanishes to infinite order, but the two arguments
of $\psi^{-1}$ vanish to finite order, we can only conclude that the lefthand side of $(\star)$ tends to some
positive constant $\gamma$.  Hence
$$
|\lambda_j(z_1, z_2, \dots, z_{n-1}| = \gamma
$$
on the boundary of the unit ball in $(n-1)$-dimensional complex space.   By our lemma, we conclude that
$\lambda$ is a constant function.   It follows that the mappings $\varphi_j$ cannot exist (as automorphisms
of $\widetilde{\Omega}$).  Hence our proof is complete.
\endpf 
\smallskip \\

\begin{remark} \rm
In fact the proposition shows that there are uncountably many biholomorphically distinct domains that
satisfy the conclusion of the assertion.  For one has great freedom in selecting the function $\psi$, and
the theory of the Chern-Moser invariants establishes the biholomorphic distinctness for the different
choices.
\end{remark}

Next we have:

\begin{proposition} \sl
Consider the domain
$$
\widetilde{\widetilde{\Omega}} = \{z = (z_1, z_2, \dots, z_n) \in \CC^n: |z_1|^{2m_1} + |z_2|^{2m_2} + \cdots + |z_{n-1}|^{2m_{n-1}} + \psi(|z_n|) < 1\} \, ,
$$
where the $m_j$ are positive integers and where $\psi$ is a real-valued, even, smooth, monotone-and-convex-on-$[0,\infty)$ 
function of a real variable with $\psi(0) = 0$ that vanishes to infinite order at 0.   Then $\widetilde{\widetilde{\Omega}}$ has
compact automorphism group.
\end{proposition}
{\bf Proof:}  If all the the $m_j$ are greater than 1 then it is easy to calculate that all of the circles
$$
\{(e^{i\theta}, 0, \dots, 0): 0 \leq \theta \leq 2\pi\} \ , \ \{(0, e^{i\theta}, 0, \dots, 0): 0 \leq \theta \leq 2\pi\} \ , \ 
      \dots \ , \
$$
$$
\{(0,0, \dots, 0, e^{i\theta}): 0 \leq \theta \leq 2\pi\}
$$
consist of weakly pseudoconvex points.  In fact points on any of those circles have $(n-1)$ weakly
pseudoconvex directions.  The boundary points with all coordinates nonvanishing are strongly
pseudoconvex.   Thus points of the latter type form a relatively open and dense set in the boundary.
One may argue that, if the automorphism group were noncompact and there were a point $X \in \widetilde{\widetilde{\Omega}}$
and a point $P \in \partial \widetilde{\widetilde{\Omega}}$ with $\varphi_j(X) \ra P$ for some automorphisms $\varphi_j$, then
at least one of these circles would be moved towards $P$ by the $\varphi_j$.   As a result, weakly
pseudoconvex boundary points would be mapped to strongly pseudoconvex boundary points.   And that is impossible.

So the only case of interest is that where at least one of the $m_j$ is equal to 1.  Say for
specificity that $m_1 = 1$.   In that case points of the form $(e^{i\theta}, 0, \dots, 0)$ have
$(n-1)$ weakly pseudoconvex directions, whereas points on the other indicated circles have one weakly
pseudoconvex direction.  As before, points with all coordinates nonvanishing are strongly
pseudoconvex.  So we can argue as above to see that the automorphism group must be compact.
\endpf 
\smallskip \\

Another type of domain that we may consider is described in the next proposition:

\begin{proposition} \sl
Consider a domain of the form
$$
\Omega' = \{(z_1, z_2, \dots, z_n) \in \CC^n: |z_1|^2 + 2 \exp(-1/(|z_2|^2 + \cdots + |z_n|^2)\} < 1 \, .
$$
Then $\Omega'$ must have compact automorphism group.
\end{proposition}
{\bf Proof:}  As usual, points with $z_2 \ne 0$,$z_3 \ne 0$, \dots, $z_n \ne 0$ are strongly pseudoconvex.  By the Bun Wong/Rosay
theorem, they cannot be orbit accumulation points.

If a point of the form $(e^{i\theta}, 0, \dots, 0)$ is an orbit accumulation point, then
we may argue as in the proof of the first proposition.  
The main difference is that line $(\dagger)$ is replaced by
$$
\left | \frac{z_1 - a_j^1}{1 - \overline{a}_j^1 z_1} \right |^2 + \rho(z_2, z_3, \dots, z_n) \cdot \lambda_j(z_1) < 1 \, ,
$$
where $\rho$ is a unitary rotation in $(n-1)$ variables.   The argument is now completed as before.
\endpf
\smallskip \\

\section{Concluding Remarks}

There is no Riemann mapping theorem in several complex variables.  Thus
we seek other means of contrasting and comparing domains in $\CC^n$.  The
automorphism group and its associated features have proved to be powerful
and flexible invariants in this study.   Certainly the Levi geometry of
boundary orbit accumulation points has been studied extensively ([GRK2] and [ISK]).
The results of this paper fit into that program.

We hope to continue these studies in future papers.

\newpage

\null \vspace*{-1in}

\noindent {\Large \sc References}
\bigskip  \\

\begin{enumerate}

\item[{\bf [BEB]}] S. Bell and H. Boas, Regularity of the
Bergman projection in weakly pseudoconvex domains, {\it Math.\
Ann.} 257(1981), 23--30.

\item[{\bf [CHM]}]  S.-S. Chern and J. Moser, Real hypersurfaces in complex
manifolds, {\em Acta Math.} 133(1974), 219-271.

\item[{\bf [GRK1]}] R. E. Greene and S. G. Krantz, Techniques
for studying automorphisms of weakly pseudoconvex domains,
{\it Several Complex Variables} (Stockholm, 1987/1988),
Math.\ Notes 38, Princeton University Press, Princeton, NJ,
1993, 389--410.

\item[{\bf [GRK2]}] R. E. Greene and S. G. Krantz,
Biholomorphic self-maps of domains, {\it Complex Analysis II}
(C. Berenstein, ed.), Springer Lecture Notes, vol. 1276, 1987,
136-207.	

\item[{\bf [HEL]}] S. Helgason, {\it Differential Geometry and
Symmetric Spaces}, Academic Press, New York, 1962.

\item[{\bf [ISK]}] A. Isaev and S. G. Krantz, Domains with
non-compact automorphism group: a survey, {\it Adv.\ Math.}
146(1999), 1--38.

\item[{\bf [KIM]}]  K.-T. Kim, Domains in $\CC^n$ with a
piecewise Levi flat boundary which possess a noncompact
automorphism group, {\it Math. Ann.} 292(1992), 575--586.

\item[{\bf [KIK]}] K.-T. Kim and S. G. Krantz, Complex scaling
and domains with non-compact automorphism group, {\it Illinois
Journal of Math.} 45(2001), 1273--1299.
			  
\item[{\bf [KRA1]}]  S. G. Krantz, {\it Function Theory of Several Complex Variables},
$2^{\rm nd}$ ed., American Mathematical Society, Providence, RI, 2001.
\item[{\bf [NAR]}] R. Narasimhan, {\it Several Complex
Variables}, University of Chicago Press, Chicago, 1971.

\item[{\bf [ROS]}] J.-P. Rosay, Sur une characterization de la
boule parmi les domains de $\CC^n$ par son groupe
d'automorphismes, {\it Ann. Inst. Four. Grenoble} XXIX(1979),
91-97.

\item[{\bf [WON]}] B. Wong, Characterizations of the ball in
$\CC^n$ by its automorphism group, {\it Invent. Math.}
41(1977), 253-257.

\end{enumerate}
\vspace*{.95in}

\small

\noindent \begin{quote}
Department of Mathematics \\
Washington University in St. Louis \\
St.\ Louis, Missouri 63130 \\ 			   
{\tt sk@math.wustl.edu}
\end{quote}

\end{document}